\documentclass[12pt,a4paper]{amsart}
\usepackage{amsfonts,color}
\usepackage{amsthm}
\usepackage{amsmath}
\usepackage{amscd}
\usepackage[latin2]{inputenc}
\usepackage{t1enc}
\usepackage[mathscr]{eucal}
\usepackage{indentfirst}
\usepackage{graphicx}
\usepackage{graphics}
\usepackage{pict2e}
\usepackage{epic}
\usepackage{url}
\usepackage{epstopdf}
\usepackage{comment}

\usepackage{tikz}

\usetikzlibrary{calc,decorations.pathreplacing}
\usetikzlibrary{quotes,angles}

\tikzset{point/.style={circle,fill,draw,inner sep=0,minimum size=1pt}}
\tikzset{vertex/.style={circle,fill,draw,inner sep=0,minimum size=7pt}}
\tikzset{overtex/.style={circle,fill=none,draw,inner sep=0,minimum size=7pt}}

\numberwithin{equation}{section}
\usepackage[margin=2.6cm]{geometry}

\allowdisplaybreaks

\theoremstyle{plain}
\newtheorem{Th}{Theorem}[section]
\newtheorem{Lemma}[Th]{Lemma}

\newtheorem{Prop}[Th]{Proposition}

 \theoremstyle{definition}
\newtheorem{Def}[Th]{Definition}
\newtheorem{Conj}[Th]{Conjecture}
\newtheorem{Rem}[Th]{Remark}
\newtheorem{?}[Th]{Problem}

\renewcommand{\P}{\mathbb{P}}

\begin{document}

\title{Upper bound for the number of spanning forests of regular graphs}

\author[F. Bencs]{Ferenc Bencs}

\address{Korteweg de Vries Institute for Mathematics, University of Amsterdam. P.O. Box 94248 1090 GE Amsterdam  The Netherlands.}

\email{ferenc.bencs@gmail.com}

\author[P. Csikv\'ari]{P\'{e}ter Csikv\'{a}ri}

\address{Alfr\'ed R\'enyi Institute of Mathematics, H-1053 Budapest, Re\'altanoda utca 13-15}

\email{peter.csikvari@gmail.com}

\thanks{The first author is supported by the NKFIH (National Research, Development and Innovation Office, Hungary) grant KKP-133921. The second author  is partially supported by the  Counting in Sparse Graphs Lend\"ulet Research Group.
}



\begin{abstract} 
We show that if $G$ is a $d$--regular graph on $n$ vertices, then the number of spanning forests $F(G)$ satisfies $F(G)\leq d^n$. The previous best bound due to Kahale and Schulman gave $(d+1/2+O(1/d))^n$. We also have the more precise conjecture that
$$F(G)^{1/n}\leq \frac{(d-1)^{d-1}}{(d^2-2d-1)^{d/2-1}}.$$
If this conjecture is true, then the expression on the right hand side is the best possible. 
\end{abstract}

\maketitle

\section{Introduction}

In this paper, we study the number of spanning forests of regular graphs. Let $\tau(G)$ and $F(G)$ denote the number of spanning trees and spanning forests, respectively. By a spanning forest $F$ we mean an acyclic set of edges. The word spanning simply refers to the fact that we always consider $V(G)$ to be the vertex set of the forest $F$.

One expects similar theorems to be true for $\tau(G)$ and $F(G)$, still, it turns out that they are very different in nature. It is known that $\tau(G)$ can be computed in polynomial time by the matrix-tree theorem \cite{kirchhoff1847ueber}. On the other hand, $F(G)$ is \#P-hard to compute \cite{jaeger1990computational}. In fact, for studying $\tau(G)$ one can use the whole machinery of Laplacian matrices and eigenvalues, whereas for $F(G)$ there is  no such toolkit. Nevertheless one expects similar or at least analogous behaviours in certain problems even if we cannot prove them for spanning forests. For instance, it is  known \cite{feder1992balanced,kirchhoff1847ueber} that if $e,f\in E(G)$ are fixed edges of a graph $G$, and $\textbf{T}$ is a random spanning tree chosen uniformly at random, then
$$\mathbb{P}(e,f\in \textbf{T})\leq \mathbb{P}(e\in \textbf{T})\mathbb{P}(f\in \textbf{T}).$$
For a uniformly chosen random spanning forest $\textbf{F}$ the similar statement
$$\mathbb{P}(e,f\in \textbf{F})\leq \mathbb{P}(e\in \textbf{F})\mathbb{P}(f\in \textbf{F})$$
is a well-known open problem \cite{grimmett2004negative,pemantle2000towards}.

A classical result due to McKay gives an upper bound on the number of spanning trees of regular graphs.

\begin{Th}[McKay \cite{mckay1983spanning}]
Let $d\geq 3$ and let $G$ be a $d$--regular graph on $n$ vertices. There exists an absolute constant $c_d$ such that
$$\tau(G)\leq \frac{c_d \ln(n)}{n}\left(\frac{(d-1)^{d-1}}{(d^2-2d)^{d/2-1}}\right)^n.$$
\end{Th}

In this paper we conjecture the following analogue of this theorem.

\begin{Conj} \label{forest-conj}
Let $d\geq 3$ and let $G$ be a $d$--regular graph on $n$ vertices. Then
$$F(G)\leq \left(\frac{(d-1)^{d-1}}{(d^2-2d-1)^{d/2-1}}\right)^n.$$
\end{Conj}

This conjecture is motivated by the following theorem of the authors.

\begin{Th}[Bencs and Csikv\'ari \cite{bencs2021evaluations}] \label{forest-limit}
Let $v(H)$ denote the number of vertices of a graph $H$. Let $F(H)$ denote the number of spanning forests of the graph $H$. Furthermore, let $g(H)$ denote the girth of the graph $H$, that is, the length of the shortest cycle. Let $(G_n)_n$ be a sequence of $d$--regular graphs, where $d\geq 3$, such that  $\lim_{n\to \infty} g(G_n)=\infty$.  Then
$$\lim_{n\to \infty}F(G_n)^{1/v(G_n)}=\frac{(d-1)^{d-1}}{(d^2-2d-1)^{d/2-1}}.$$
If $(G_n)_n$ is a sequence of random $d$--regular graphs with growing number of vertices, then the same statement holds true asymptotically almost surely.
\end{Th}

Although we are not able to prove Conjecture~\ref{forest-conj}, we can prove the following result.

\begin{Th} \label{forest-theorem} Let $d\geq 2$ be an integer.
Let $G$ be a $d$--regular graph on $n$ vertices. Then
$F(G)\leq d^n$.
\end{Th}

To compare Conjecture~\ref{forest-conj} and Theorem~\ref{forest-theorem} we note that
$$\frac{(d-1)^{d-1}}{(d^2-2d-1)^{d/2-1}}=d-\frac{1}{2d}-\frac{1}{3d^2}-\frac{1}{8d^3}+O\left(\frac{1}{d^4}\right). $$
We also note that the previous best bound on $F(G)$ is due to Kahale and Schulman \cite{kahale1996bounds}, and they proved that $F(G)\leq C_d^n$, where
$$C_d=d+\frac{1}{2}+\frac{1}{8d}+\frac{13}{48d^2}+O\left(\frac{1}{d^3}\right).$$
Thus, the actual improvement in the constant is more than $1/2$. Note that Borb\'enyi, Csikv\'ari and Luo \cite{borbenyi2020number} were able to improve on the results of Kahale and Schulman for certain values of $d$, but even those bounds are far from the result of Theorem~\ref{forest-theorem}. Let us also mention that we will actually prove a slightly stronger result, namely,
$$F(G)^{1/v(G)}\leq \sqrt{2}\mu_{K_{d+1}}\left(\frac{d+1}{\sqrt{2}}\right)^{1/(d+1)},$$
where $\mu_G(x)$ is the matching polynomial of a graph $G$ (see Section~\ref{sec: matching polynomial}), and $K_{d+1}$ is the complete graph on $d+1$ vertices. A considerable work is actually needed to prove that the right hand side is at most $d$ if $d\geq 4$.
\bigskip

One may wish to compare our results to various specific graph classes. One notable special case concerns  the number of forests of Archimedean lattices. In this case we recommend the short survey of Chang and Shrock \cite{chang2020exponential}.
\bigskip

\noindent \textbf{Notation.} Throughout the paper $G=(V,E)$ denotes a graph with vertex set $V$ and edge set $E$. We will denote $|E(G)|$ by $e(G)$, and $|V(G)|$ by $v(G)$. If we fix a graph $G$, then we use $v(G)=n$. The degree of a vertex $v$ is the number of neighbors of $v$, and is denoted by $d_v$. 
$K_n$ denotes the complete graph on $n$ vertices.
\bigskip

\noindent \textbf{Organization of this paper.} In the next section we study a polynomial that we call $R_G(z)$. This polynomial was introduced in  paper \cite{bencs2021evaluations}. It turns out that $F(G)\leq R_G(2)$. Since $R_G(z)$ is very strongly connected to the so-called matching polynomial we also collect all relevant results about the matching polynomial here. In Section~\ref{sec: proof of main theorem} we prove  Theorem~\ref{forest-theorem}. In Section~\ref{sec: conjecture revisited} we revisit Conjecture~\ref{forest-conj}.

\section{The polynomial $R_G(z)$} \label{sec: polynomial $R_G(z)$}

In this section we study a polynomial that is related to the matching polynomial of a graph $G$. Recall that a set of edges $M$ is a matching if  no two edges of $M$ have a common end vertex. A $k$-matching is simply a matching of size $k$.

\begin{Def}
Let
$$R_G(z)=\sum_{M\in \mathcal{M}(G)}(-z)^{|M|}\prod_{v\notin V(M)}(z+d_v-1),$$
where $d_v$ is the degree of the vertex $v$, and $\mathcal{M}(G)$ is the set of matchings of $G$, including the empty one.
\end{Def}

\begin{Rem}
This polynomial was introduced by the authors in  paper \cite{bencs2021evaluations} in order to study the Tutte polynomial. In this paper we strongly build on \cite{bencs2021evaluations}.
\end{Rem}

Next we need the concept of pseudo-forests for a better understanding of $R_G(z)$.

\begin{Def}
We say that an $A\subseteq E(G)$ is a pseudo-forest of $G$ if each of its connected components are forests or unicyclic graphs. Let $\mathcal{PF}(G)$ be the set of pseudo-forests of $G$. For a pseudo-forest $A$ let $c(A)$ be the number of cycles in $A$. \end{Def}

The main result about $R_G(z)$ is the following characterization in terms of pseudo-forests.

\begin{Th}[Bencs and Csikv\'ari \cite{bencs2021evaluations}] \label{R-pseudo-forests}
For a graph $G$ we have
$$R_G(z+1)=\sum_{k=0}^n\left(\sum_{A\in \mathcal{PF}(G) \atop |A|=k}2^{c(A)}\right)z^{n-k}.$$
\end{Th}

Let us introduce the polynomial
$$F_G(z)=\sum_{k=0}^nf_k(G)z^{n-k},$$
where $f_k(G)$ denotes the number of forests with exactly $k$ edges. It turns out that $F_G(z)$ is a specialization of the  
Tutte polynomial \cite{tutte1954contribution}. The Tutte polynomial  of a graph $G$ is defined as follows:
$$T_G(x,y)=\sum_{A\subseteq E(G)}(x-1)^{k(A)-k(G)}(y-1)^{k(A)+|A|-|V|},$$
where $k(A)$ is the number of connected components of the graph with vertex set $V=V(G)$ and edge set $A$, and $k(G)$ is the number of connected components of the graph $G$. For a comprehensive introduction to the Tutte polynomial see the surveys \cite{brylawski1992tutte,crapo1969tutte,ellis2011graph,welsh1999tutte}.   

It is easy to see that
$$F_G(z)=z^{k(G)}T_G(z+1,1).$$
Since a forest is a pseudo-forest without any cycle we immediately have
$$F_G(z)\leq R_G(z+1)$$
for positive $z$. In particular, 
\begin{equation}\label{f_versus_r}
F(G)=F_G(1)\leq R_G(2)    
\end{equation}

\subsection{Matching polynomial} \label{sec: matching polynomial}

In this section, we collect a few facts about the so-called matching polynomial. The matching polynomial of a graph $G$ is defined as follows. Let 
$$\mu_G(z)=\sum_{k=0}^{n/2}(-1)^km_k(G)z^{n-2k},$$
where $m_k(G)$ denotes the number of matchings of size $k$.
Note that for a $d$--regular graph $G$ we have
\begin{equation}\label{r_as_matching}
R_G(z)=\sum_{M\in \mathcal{M}(G)}(-z)^{|M|}(d+z-1)^{n-2|M|}=z^{n/2}\mu_G\left(\frac{d-1+z}{\sqrt{z}}\right).
\end{equation}
The following theorem about the matching polynomial is fundamental for us.

\begin{Th}[Heilmann and Lieb \cite{heilmann1972theory}] \label{th: HeiLie}
All zeros of the matching polynomial $\mu_G(z)$ are real. Furthermore, if the largest degree $\Delta$ satisfies $\Delta \geq 2$, then all zeros lie in the interval $(-2\sqrt{\Delta-1},2\sqrt{\Delta-1})$.
\end{Th}

Let $\mu_G(z)=\prod_{i=1}^n(z-\alpha_i)$, and
$s_k(G)=\sum_{i=1}^n \alpha_i^k$.
The quantities $s_k(G)$ have a combinatorial meaning. They count the so-called tree-like walks. Before introducing the concept of tree-like walks, we first require the definition of a path-tree.

\begin{figure} 
\begin{center}
\includegraphics[scale=0.8]{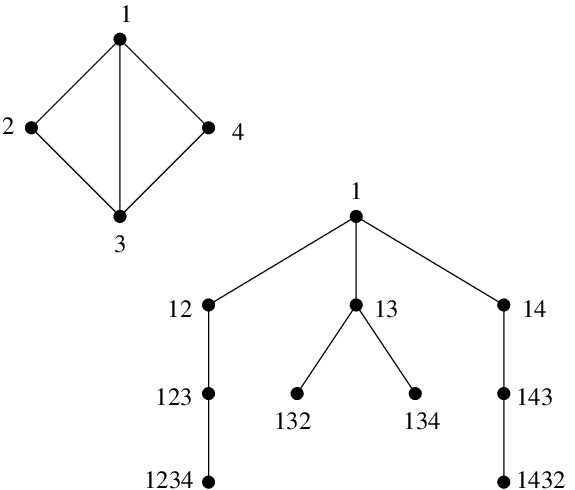}
\caption{A path-tree from vertex $1$.}  \label{fig: path-tree}
\end{center}
\end{figure}

\begin{Def} Let $G$ be graph with a given vertex $u$. The \textit{path-tree} $T(G,u)$ is defined as follows. The vertices of $T(G,u)$ are the paths in $G$ which start at the  vertex $u$ and two paths as vertices are joined by an edge if one of them is a one-step extension of the other. See Figure~\ref{fig: path-tree}.
\end{Def}

\begin{Def}
A closed tree-like walk of length $\ell$ is a closed walk on $T(G,u)$ of length $\ell$ starting and ending at the vertex $u$.
\end{Def}

Note that a priori the tree-like walk is on the path-tree, although one can make a correspondence with certain walks on the graph itself. Indeed, a walk in the tree $T(G,u)$ from $u$ can be imagined as follows. Suppose that in the graph $G$ a worm is sitting at  vertex $u$ at the beginning. Then at each step the worm can either grow or pull back its head. When it grows it can move its head to a neighboring unoccupied vertex while keeping its tail at vertex $u$.  At each step the worm occupies a path in the graph $G$. A closed walk in the tree $T(G,u)$ from $u$ to $u$ corresponds to the case when at the final step the worm occupies only vertex $u$.  Let $s_k(G,u)$ be the number of tree-like closed walks of length $k$ starting and ending at vertex $u$. Note that for different vertices $u$ and $v$ the quantities $s_k(G,u)$ and $s_k(G,v)$ may be different.

\begin{Lemma}[Godsil \cite{godsil1993algebraic}] \label{tree-like walks}
The number of all closed tree-like walks of length $\ell$ is exactly $s_{\ell}(G)$, that is, $s_{\ell}(G)=\sum_{u\in V(G)}s_{\ell}(G,u)$. 
\end{Lemma}

The following inequalities will be useful when we prove Theorem~\ref{forest-theorem}.

\begin{Lemma} \label{comparison2} For every $d$--regular graph $G$ and $k\geq 1$ we have 
$$\frac{s_k(G)}{v(G)}\geq \frac{s_k(K_{d+1})}{v(K_{d+1})}.$$
\end{Lemma}

\begin{proof}
The path-tree of the complete graph on $d+1$ vertices is a subtree of the path-tree of any regular graph. Indeed, if we have  a path $P$ in $G$ consisting of $t$ vertices, then we can extend it to at least $d+1-t$ ways to a path on $t+1$ vertices. In $K_{d+1}$ we have equality here. 

From this it follows that if $u\in V(G)$, then $s_k(G,u)\geq s_k(K_{d+1},r)$, where $r$ is any vertex of $K_{d+1}$. Then
$$\frac{s_k(G)}{v(G)}=\frac{1}{v(G)}\sum_{u\in V(G)}s_k(G,u)\geq s_{k}(K_{d+1},r)=\frac{s_k(K_{d+1})}{v(K_{d+1})}.$$
\end{proof}

\begin{Lemma} \label{integral1} Let $G$ be a $d$--regular graph with $v(G)$ vertices, and let $a>2\sqrt{d-1}$. Then we have
$$\frac{1}{v(G)}\log \mu_G(a)\leq \frac{1}{v(K_{d+1})}\log \mu_{K_{d+1}}(a).$$
\end{Lemma}

\begin{proof} For a graph $G$ we can write the matching polynomial as $\mu_G(z)=\prod_{i=1}^{v(G)}(z-\alpha_i)$, where $|\alpha_i|<2\sqrt{d-1}<a$ by Theorem~\ref{th: HeiLie}. Then we have
$$\frac{1}{v(G)}\log \mu_G(a)=\frac{1}{v(G)}\sum_{i=1}^{v(G)}\log(a-\alpha_i)=\log(a)+\frac{1}{v(G)}\sum_{i=1}^{v(G)}\log\left(1-\frac{\alpha_i}{a}\right)=$$
$$=\log(a)-\frac{1}{v(G)}\sum_{k=1}^{\infty}\frac{1}{ka^k}\sum_{i=1}^n\alpha_i^k=\log(a)-\sum_{k=1}^{\infty}\frac{1}{ka^k}\frac{s_k(G)}{v(G)}.$$
The upper bound now follows using that
$\frac{s_k(G)}{v(G)}\geq \frac{s_k(K_{d+1})}{v(K_{d+1})}$.
\end{proof}

\section{Proof of Theorem~\ref{forest-theorem}} \label{sec: proof of main theorem}

In this section we prove Theorem~\ref{forest-theorem}. Note that the cases $d=2$ and $d=3$ are trivial since for $d=2$
$$F(G)\leq 2^{e(G)}=2^{dn/2}=2^n,$$
and for $d=3$ we have
$$F(G)\leq 2^{e(G)}=2^{3n/2}=(2\sqrt{2})^n<3^n.$$
From now on we can assume that $d\geq 4$.
By inequality~\ref{f_versus_r} and identity~\ref{r_as_matching} we have  that if $G$ is a $d$--regular graph on $n$ vertices, then
$$F(G)\leq R_G(2)=2^{n/2}\mu_G\left(\frac{d+1}{\sqrt{2}}\right).$$
Hence it is enough to prove that if $d\geq 4$, then
$$2^{n/2}\mu_G\left(\frac{d+1}{\sqrt{2}}\right)<d^n.$$
First we show that the quantity
$\sqrt{2}\mu_G\left(\frac{d+1}{\sqrt{2}}\right)^{1/{v(G)}}$
is maximised by the graph $K_{d+1}$ among $d$--regular graphs. Then we show that
$$\sqrt{2}\mu_{K_{d+1}}\left(\frac{d+1}{\sqrt{2}}\right)^{1/(d+1)}<d$$
if $d\geq 4$.
Let us start with the first inequality.

\begin{Prop}
Let $d\geq 4$. For a $d$--regular graph $G$ on $n$ vertices we have
$$\sqrt{2}\mu_G\left(\frac{d+1}{\sqrt{2}}\right)^{1/n}\leq \sqrt{2}\mu_{K_{d+1}}\left(\frac{d+1}{\sqrt{2}}\right)^{1/(d+1)}.$$

\end{Prop}

\begin{proof}
This follows from Lemma~\ref{integral1} since $\frac{d+1}{\sqrt{2}}> 2\sqrt{d-1}$ for every $d\geq 4$.
\end{proof}

So we have 
$$F(G)<\left(\sqrt{2}\mu_{K_{d+1}}\left(\frac{d+1}{\sqrt{2}}\right)^{1/(d+1)}\right)^n$$
for every $d$--regular graph $G$ on $n$ vertices.
\bigskip

\begin{table}
\begin{center}
\begin{tabular}{|c|c|c|c|} \hline
$d$ & $\frac{(d-1)^{d-1}}{(d^2-2d-1)^{d/2-1}}$ & $\sqrt{2}\cdot \mu_{K_{d+1}}\left(\frac{d+1}{\sqrt{2}}\right)^{1/(d+1)}$ & $d-\frac{1}{2d}$ \\ \hline
4 & 3.85714285714286 & 3.91947904192452 & 3.87500000000000\\ \hline
5 & 4.88706270925576 & 4.91723641784679 & 4.90000000000000\\ \hline
6 & 5.90737240075614 & 5.92330714974640 & 5.91666666666667\\ \hline
7 & 6.92165915952326 & 6.93081266948571 & 6.92857142857143\\ \hline
8 & 7.93218265702205 & 7.93782732697178 & 7.93750000000000\\ \hline
9 & 8.94023598867791 & 8.94392512794576 & 8.94444444444444\\ \hline
10 & 9.94659007255980 & 9.94911825176517 & 9.95000000000000\\ \hline
11 & 10.9517282616543 & 10.9535296122445 & 10.9545454545455\\ \hline
12 & 11.9559676088814 & 11.9572931718296 & 11.9583333333333\\ \hline
13 & 12.9595242435679 & 12.9605264078823 & 12.9615384615385\\ \hline
14 & 13.9625503782448 & 13.9633255787092 & 13.9642857142857\\ \hline
15 & 14.9651562263964 & 14.9657677145412 & 14.9666666666667\\ \hline
16 & 15.9674235136022 & 15.9679140832706 & 15.9687500000000\\ \hline
17 & 16.9694141030687 & 16.9698135007115 & 16.9705882352941\\ \hline
18 & 17.9711756732928 & 17.9715050723952 & 17.9722222222222\\ \hline
19 & 18.9727455571083 & 18.9730203479578 & 18.9736842105263\\ \hline
20 & 19.9741533998487 & 19.9743849792194 & 19.9750000000000\\ \hline
\end{tabular}
\end{center}
\caption{Numerical values of the constants appearing in this paper for $4\le d\le 20$.}
\end{table}

As Table 1 suggests this  bound is actually better than $d-\frac{1}{2d}$ from $d\geq 9$, and for all $d$ it is less than $d$. We only prove this latter result.

\begin{Lemma}
For $d\geq 4$ we have
$$\sqrt{2}\mu_{K_{d+1}}\left(\frac{d+1}{\sqrt{2}}\right)^{1/(d+1)}<d.$$
\end{Lemma}

The following proof is quite long and technical, but completely standard analysis.

\begin{proof}
It will be convenient to introduce the notation $n=d+1$. Then the required inequality is equivalent with
$$2^{n/2}\mu_{K_n}\left(\frac{n}{\sqrt{2}}\right)<(n-1)^n$$
for $n\geq 5$. The plan is the following: we will show that if we divide by $n^n$, then
$$\frac{2^{n/2}}{n^n}\mu_{K_n}\left(\frac{n}{\sqrt{2}}\right)=\frac{1}{e}-\frac{1}{en}+O\left(\frac{1}{n^2}\right),$$
and
$$\frac{(n-1)^n}{n^n}=\frac{1}{e}-\frac{1}{2en}+O\left(\frac{1}{n^2}\right),$$
where $e=2.71..$ is the base of the natural logarithm. We will prove these statements with quantitative error bounds to determine those numbers $n$ for which the error term is less than the second main term. For smaller $n$ one may check the statement with computer.

Note that $m_k(K_n)=\frac{n!}{2^kk!(n-2k)!}$, and so
\begin{align*}
\frac{2^{n/2}}{n^n}\mu_{K_n}\left(\frac{n}{\sqrt{2}}\right)&=\frac{2^{n/2}}{n^n}\sum_{k=0}^{\lfloor n/2 \rfloor}(-1)^km_k(K_n)\left(\frac{n}{\sqrt{2}}\right)^{n-2k}\\
&=\frac{2^{n/2}}{n^n}\sum_{k=0}^{\lfloor n/2 \rfloor}(-1)^k\frac{n!}{2^kk!(n-2k)!}\left(\frac{n}{\sqrt{2}}\right)^{n-2k}\\
&=\sum_{k=0}^{\lfloor n/2 \rfloor}\frac{(-1)^k}{k!}\frac{n(n-1)\dots (n-2k+1)}{n^{2k}}\\
&=\sum_{k=0}^{\lfloor n/2 \rfloor}\frac{(-1)^k}{k!}\prod_{i=0}^{2k-1}\left(1-\frac{i}{n}\right).
\end{align*}
Similarly,
$$
\frac{(n-1)^n}{n^n}=\left(1-\frac{1}{n}\right)^n
=\sum_{k=0}^n\binom{n}{k}(-1)^k\frac{1}{n^k}
=\sum_{k=0}^n\frac{(-1)^k}{k!}\prod_{i=0}^{k-1}\left(1-\frac{i}{n}\right).
$$
In what follows we will often use the following two observations. If $b_1> b_2>\dots >b_r>0$ and $r>m$, then
$$\left|\sum_{j=m+1}^r(-1)^jb_j\right|\leq b_{m+1}.$$
We also need that if $a_1,\dots ,a_r>0$ such that $\sum_{i=1}^r a_i<1$, then
$$1-\sum_{i=1}^r a_i\leq \prod_{i=1}^r(1-a_i)<1-\sum_{i=1}^r a_i+\frac{1}{2}\left(\sum_{i=1}^r a_i\right)^2.$$
Here the lower bound can be proved by induction on $r$, while the upper bound follows from
$$\prod_{i=1}^r(1-a_i)<\prod_{i=1}^re^{-a_i}=\exp\left(-\sum_{i=1}^r a_i\right)<1-\sum_{i=1}^r a_i+\frac{1}{2}\left(\sum_{i=1}^r a_i\right)^2.$$
since $e^{-x}=1-x+\frac{x^2}{2}-\frac{x^3}{3!}+\dots$ and $\sum_{j=2}^{\infty}\left(-\frac{x^{2j-1}}{(2j-1)!}+\frac{x^{2j}}{(2j)!}\right)<0$ if $x<1$.
The sequences $b_k=\frac{1}{k!}\prod_{i=0}^{2k-1}\left(1-\frac{i}{n}\right)$ and $b'_k=\frac{1}{k!}\prod_{i=0}^{k-1}\left(1-\frac{i}{n}\right)$ are monotone decreasing, so for $m< \lfloor n/2 \rfloor$ we have
$$\sum_{k=0}^{\lfloor n/2 \rfloor}\frac{(-1)^k}{k!}\prod_{i=0}^{2k-1}\left(1-\frac{i}{n}\right)=\sum_{k=0}^{m}\frac{(-1)^k}{k!}\prod_{i=0}^{2k-1}\left(1-\frac{i}{n}\right)+R^{(1)}_m$$
and
$$\sum_{k=0}^n\frac{(-1)^k}{k!}\prod_{i=0}^{k-1}\left(1-\frac{i}{n}\right)=\sum_{k=0}^m\frac{(-1)^k}{k!}\prod_{i=0}^{k-1}\left(1-\frac{i}{n}\right)+r^{(1)}_m,$$
where $|R^{(1)}_m|,|r^{(1)}_m|<\frac{1}{(m+1)!}$. We will choose $m$ in such a way that $\sum_{j=0}^{2m-1}\frac{j}{n}=\frac{\binom{2m}{2}}{n}<1$. Then for $k\leq m$ we have
$$\prod_{i=0}^{2k-1}\left(1-\frac{i}{n}\right)=1-\frac{\binom{2k}{2}}{n}+R_k\ \ \ \mbox{and}\ \ \ 
\prod_{i=0}^{k-1}\left(1-\frac{i}{n}\right)=1-\frac{\binom{k}{2}}{n}+r_k,$$
where $|R_k|\leq \frac{1}{2}\left(\frac{\binom{2k}{2}}{n}\right)^2$ and
$|r_k|\leq \frac{1}{2}\left(\frac{\binom{k}{2}}{n}\right)^2$.
Hence
\begin{align*}
\sum_{k=0}^{\lfloor n/2 \rfloor}\frac{(-1)^k}{k!}\prod_{i=0}^{2k-1}\left(1-\frac{i}{n}\right)&=\sum_{k=0}^{m}\frac{(-1)^k}{k!}\prod_{i=0}^{2k-1}\left(1-\frac{i}{n}\right)+R^{(1)}_m\\
&=\sum_{k=0}^{m}\frac{(-1)^k}{k!}\left(1-\frac{\binom{2k}{2}}{n}+R_k\right)+R^{(1)}_m\\
&=\sum_{k=0}^{m}\frac{(-1)^k}{k!}+\frac{1}{n}\sum_{k=0}^{m}(-1)^k\frac{\binom{2k}{2}}{k!}+\sum_{k=0}^{m}\frac{(-1)^k}{k!}R_k+R^{(1)}_m\\
&=\sum_{k=0}^{m}\frac{(-1)^k}{k!}+\frac{1}{n}\left(-\frac{1}{e}+R^{(2)}_m\right)+\sum_{k=0}^{m}\frac{(-1)^k}{k!}R_k+R^{(1)}_m,
\end{align*}
where $|R^{(2)}_m|\leq \frac{\binom{2(m+1)}{2}}{(m+1)!}$.
Similarly,
\begin{align*}
\sum_{k=0}^{n}\frac{(-1)^k}{k!}\prod_{i=0}^{k-1}\left(1-\frac{i}{n}\right)&=\sum_{k=0}^{m}\frac{(-1)^k}{k!}\prod_{i=0}^{k-1}\left(1-\frac{i}{n}\right)+r^{(1)}_m\\
&=\sum_{k=0}^{m}\frac{(-1)^k}{k!}\left(1-\frac{\binom{k}{2}}{n}+r_k\right)+r^{(1)}_m\\
&=\sum_{k=0}^{m}\frac{(-1)^k}{k!}+\frac{1}{n}\sum_{k=0}^{m}(-1)^k\frac{\binom{k}{2}}{k!}+\sum_{k=0}^{m}\frac{(-1)^k}{k!}r_k+r^{(1)}_m\\
&=\sum_{k=0}^{m}\frac{(-1)^k}{k!}+\frac{1}{n}\left(-\frac{1}{2e}+r^{(2)}_m\right)+\sum_{k=0}^{m}\frac{(-1)^k}{k!}r_k+r^{(1)}_m,
\end{align*}
where $|r^{(2)}_m|\leq \frac{\binom{m+1}{2}}{(m+1)!}$.
Note that
$$
\left|\left(\frac{R^{(2)}_m}{n}+\sum_{k=0}^{m}\frac{(-1)^k}{k!}R_k+R^{(1)}_m\right)-\left(\frac{r^{(2)}_m}{n}+\sum_{k=0}^{m}\frac{(-1)^k}{k!}r_k+r^{(1)}_m\right)\right|$$
$$\leq \frac{|R^{(2)}_m|}{n}+\sum_{k=0}^{m}\frac{|R_k|}{k!}+|R^{(1)}_m|+\frac{|r^{(2)}_m|}{n}+\sum_{k=0}^{m}\frac{r_k}{k!}+|r^{(1)}_m|$$
$$\leq \frac{\binom{2(m+1)}{2}}{n\cdot (m+1)!}+\frac{1}{2n^2}\sum_{k=0}^{m}\frac{\binom{2k}{2}^2}{k!}+\frac{1}{(m+1)!}+\frac{\binom{m+1}{2}}{n\cdot (m+1)!}+\frac{1}{2n^2}\sum_{k=0}^{m}\frac{\binom{k}{2}^2}{k!}+\frac{1}{(m+1)!}$$
$$\leq \frac{5/2\cdot m^2+7/2\cdot m+1}{n\cdot (m+1)!}+\frac{175e}{8n^2}+\frac{2}{(m+1)!}$$
since 
$$\sum_{k=0}^{\infty}\frac{\binom{2k}{2}^2}{k!}=42e\ \ \ \text{and}\ \ \ \sum_{k=0}^{m}\frac{\binom{k}{2}^2}{k!}=\frac{7e}{4}.$$
In general, if we have a polynomial $P(x)=\sum_{r=0}^t a_rx(x-1)\dots (x-r+1)$, then
$$\sum_{k=0}^{\infty}\frac{P(k)}{k!}=\sum_{r=0}^ta_r\sum_{k=r}^{\infty}\frac{1}{(k-r)!}=e\sum_{r=0}^t a_r,$$
and in our case
$$\binom{2k}{2}^2=4k(k-1)(k-2)(k-3)+20k(k-1)(k-2)+17k(k-1)+k,$$
and
$$\binom{k}{2}^2=\frac{1}{4}k(k-1)(k-2)(k-3)+k(k-1)(k-2)+\frac{1}{2}k(k-1).$$
Hence
$$\sum_{k=0}^{n}\frac{(-1)^k}{k!}\prod_{i=0}^{k-1}\left(1-\frac{i}{n}\right)-\sum_{k=0}^{\lfloor n/2 \rfloor}\frac{(-1)^k}{k!}\prod_{i=0}^{2k-1}\left(1-\frac{i}{n}\right)$$ 
$$\geq \frac{1}{2e}\cdot \frac{1}{n}-\left(\frac{5/2\cdot m^2+7/2\cdot m+1}{n\cdot (m+1)!}+\frac{175e}{8n^2}+\frac{2}{(m+1)!}\right).$$
Let $m=\lfloor \sqrt{n/2}\rfloor$, then $\binom{2m}{2}<2m^2\leq n$. If $n\geq 392$, then $m\geq 14$, and we have
$$\frac{5/2\cdot m^2+7/2\cdot m+1}{n\cdot (m+1)!}\leq \frac{7m^2}{(m+1)!\cdot n}< \frac{7}{13!\cdot n}<\frac{1}{1000n},$$
and
$$\frac{175e}{8n^2}\leq \frac{175e}{8\cdot 392\cdot n}<\frac{1}{6n},$$
and
$$\frac{2}{(m+1)!}\leq \frac{2}{13!\cdot m\cdot (m+1)}\leq \frac{2}{6000\cdot n/3}=\frac{1}{1000n}.$$
Hence
$$\frac{1}{2e}\cdot \frac{1}{n}-\left(\frac{5/2\cdot m^2+7/2\cdot m+1}{n\cdot (m+1)!}+\frac{7e}{4n^2}+\frac{2}{(m+1)!}\right)\geq \left(\frac{1}{2e}-\frac{1}{1000}-\frac{1}{6}-\frac{1}{1000}\right)\frac{1}{n}>0.$$
So for $n\geq 392$ we have 
$$\frac{2^{n/2}}{n^n}\sum_{k=0}^{\lfloor n/2 \rfloor}(-1)^km_k(K_n)\left(\frac{n}{\sqrt{2}}\right)^{n-2k}<\frac{(n-1)^n}{n^n}.$$
The authors checked  the cases $5\leq n\leq 392$ with a computer. For courtesy we made available the Sagemath code of the computation at \cite{bencs2022Sagemath}.
\end{proof}

\section{Conjecture~\ref{forest-conj} revisited} \label{sec: conjecture revisited}

In this section we revisit Conjecture~\ref{forest-conj}. Here we review a result of Borb\'enyi, Csikv\'ari and Luo \cite{borbenyi2020number}. They showed that a well-known conjecture mentioned already in the introduction actually implies Conjecture~\ref{forest-conj}. This well-known conjecture is about a negative correlation inequality. The precise statement is the following.

\begin{Conj}[\cite{grimmett2004negative, pemantle2000towards}]\label{correlation-conjecture} Let $G$ be a graph and let $\textbf{F}$ be a random forest chosen uniformly from all the forests of $G$. Let $e,f\in E(G)$, then
$$\P(e,f\in \textbf{F})\leq \P(e\in \textbf{F})\P(f\in \textbf{F}).$$
\end{Conj}

Assuming Conjecture~\ref{correlation-conjecture}  Borb\'enyi, Csikv\'ari and Luo proved a result on forests of $2$-covers which then implies Conjecture~\ref{forest-conj}.

Recall that a $k$-cover (or $k$-lift) $H$ of a graph $G$ is defined as follows. The vertex set of  $H$ is $V(H)=V(G)\times \{0,1,\dots, k-1\}$, and the edge set of $H$ is $E(H)=\cup_{(u,v)\in E(G)}M_{u,v}$, where $M_{u,v}$ is a perfect matching between the vertex set $L_u=\{(u,i)\ |\ 0\leq i\leq k-1\}$ and $L_v=\{(v,i)\ |\ 0\leq i\leq k-1\}$.

\begin{figure}[h!]
\begin{center}
\scalebox{.8}{\includegraphics{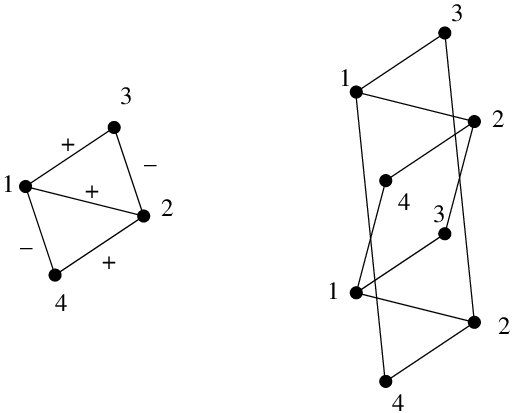}}
\caption{A $2$-lift.}
\end{center}
\end{figure}

When $k=2$ one can encode a $2$-lift $H$ by putting signs on the edges of the graph $G$: the $+$ sign means that we use the matching $((u,0),(v,0)),((u,1),(v,1))$ at the edge $(u,v)$,
the $-$ sign means that we use the matching  $((u,0),(v,1)),((u,1),(v,0))$ at the edge $(u,v)$. For instance, if we put $+$ signs to every edge, then we simply get the disjoint union $G\cup G$, and if we put $-$ signs everywhere, then the obtained $2$-cover $H$ is simply the tensor product $G\times K_2$. Observe that if $G$ is bipartite, then $G\cup G=G\times K_2$, but other $2$-covers might differ from $G\cup G$.

\begin{Th}[Borb\'enyi, Csikv\'ari and Luo \cite{borbenyi2020number}] \label{forest-cover}
Let $G$ be a graph, and let $H$ be a $2$-cover of $G$. If Conjecture~\ref{correlation-conjecture} is true, then we have
$$F(G\cup G)\leq F(H).$$
Thus, $F(G)^{1/v(G)}\leq F(H)^{1/v(H)}$.
\end{Th}

There is a nice property of covers that is related to the girth. For every graph $G$, there is a sequence of graphs $(G_n)_n$ such that $G_0=G$, $G_k$ is a $2$-cover of $G_{k-1}$, and  $g(G_k)\to \infty$. This is an observation due to Linial \cite{linial}, his proof is also given in \cite{csikvari2017lower}. This observation, Theorem~\ref{forest-limit} and Theorem~\ref{forest-cover} (assuming Conjecture~\ref{correlation-conjecture}) together implies that 
$$\sup_{G\in \mathcal{G}_d}F(G)^{1/v(G)}=\frac{(d-1)^{d-1}}{(d^2-2d-1)^{d/2-1}},$$
where $\mathcal{G}_d$ denotes the family of $d$-regular graphs.
\bigskip

One might wonder whether an analogue of Theorem~\ref{forest-theorem} is true for non-regular graphs too. The most natural guess 
$$F(G)\leq \prod_{v\in V(G)}d_v$$
can only be true if we require some minimum degree condition. We speculate that this inequality can be true if the minimum degree is at least $3$. Minimum degree $2$ is not enough as gluing together $r$ cycles of length $k$ at a fixed vertex shows that for the obtained graph $G$ we have
$F(G)>\prod_{v\in V(G)}d_v$ if $(2^k-1)^r>(2r)2^{r(k-1)}$ which is true for every fixed $k$ and sufficiently large $r$. 
Note that the inequality $F(G)\leq \prod_{v\in V(G)}(d_v+1)$ is easily seen to be true for every graph $G$ (cf. \cite{thomassen2010spanning}), but gives an exponentially weaker bound for sparse graphs.
\bigskip

\noindent \textbf{Acknowledgment.} The authors are very grateful to the reviewers for the comments and suggestions, especially for pointing out a computational error in the original manuscript.

\bibliography{forest_reference}
\bibliographystyle{plain}

\end{document}